\documentclass[review]{elsarticle}

\usepackage{lineno,hyperref}
\modulolinenumbers[5]

\usepackage{url}
\usepackage{kbordermatrix}
\usepackage{amsmath,amssymb}
\usepackage{amsthm,thmtools}
\usepackage{amsfonts} 

\usepackage[dvips]{epsfig}      % ******* for figures ***************
\usepackage{graphicx}
\usepackage{enumerate}
\usepackage{color}                  %******** for colored remarks********

\usepackage{verbatim}

\usepackage{soul} % add strikeout fonts

\usepackage{amssymb}
\usepackage{amsbsy}
\usepackage{amsmath}
\usepackage{setspace}
\usepackage{url}

\usepackage{subcaption}
\usepackage{float}

\usepackage{longtable}

 \newcommand{\diag}{\operatorname{diag}}

  \newcommand{\Trace}{\operatorname{trace}}

\renewcommand{\kbldelim}{(}% Left delimiter
\renewcommand{\kbrdelim}{)}% Right delimiter

\declaretheorem[name={Example},qed={\lower-0.3ex\hbox{$\square$}} ] {Example}

\declaretheorem[name={Definition}  ] {Definition}
\declaretheorem[name={Theorem}  ] {Theorem}

\declaretheorem[name={Remark}  ] {Remark}
\declaretheorem[name={Corollary}  ] {Corollary}

\declaretheorem[name={Proposition}  ] {Proposition}

\newcommand {\R}{\mathbb R}

\newcommand{\be}{\begin{equation}}
\newcommand{\ee}{\end{equation}}

\newcommand{\myc}[1]{{\bf{#1}}}  % use this to add visible comments, by simply typing \myc{<text>}

\setstcolor{blue} % the strikeout line color 

\newcommand{\updt}[1]{{\color{black}#1}}

\usepackage{lineno}
\begin{document}
%add this line after \begin{document} but before \titlepage

\begin{frontmatter}

\title{On the Exponent of Several Classes of Oscillatory~Matrices\tnoteref{t1}}

 \author[main1]{Yoram Zarai}% and  Michael Margaliot}
%\author[main1]{Yoram Zarai and Michael Margaliot\corref{correspauth}}% and  Michael Margaliot}
\author[main1]{Michael Margaliot\corref{correspauth}}

\address[main1]{School of Elec. Eng.-Systems, Tel-Aviv University, Tel-Aviv~69978, Israel.}
\cortext[correspauth]{Corresponding author.}
\ead{michaelm@tauex.tau.ac.il}

\tnotetext[t1]{Research  supported in part by  research grants from  the Israel Science Foundation and the US-Israel Binational Science Foundation.}

%%%%%\maketitle
%%%%%%%%%%%%%%%%%%%%%%%%%%%%%%%%%%%%%%%%%%%%%%%%%%%%%%%%%%

 % AMS subject classifications (used in AMS journals)

\begin{abstract}
%%%%%%%%%%%%%%%%%
Oscillatory matrices were introduced  in the seminal work of Gantmacher and~Krein. 
An~$n\times n$
 matrix~$A$ is called oscillatory if all its minors are nonnegative and there exists a positive integer~$k$
such that all minors of~$A^k$ are positive. The smallest~$k$ for which this holds is called the exponent of the oscillatory matrix~$A$.  Gantmacher and Krein showed that the exponent is always smaller than or equal to~$n-1$. 
An important and nontrivial problem is to determine the exact value of the exponent. Here we use
the successive elementary bidiagonal factorization of oscillatory matrices, and its graph-theoretic representation, to derive an explicit expression for the exponent of several classes of oscillatory matrices, and a nontrivial  upper-bound on the exponent 
for several other classes.
\end{abstract}
%%%%%%%%%%%%%%%%%%%%%%%%%%%%%%%%%%%%%%%%%%%%%%%%%%%%%%%%%%%%%%%%%%%%%%%%%%%%

 
\end{frontmatter}


\section{Introduction}
% %%%%%%%%%%%%%
A matrix is called totally positive~(TP) [totally nonnegative~(TN)] if all its minors are positive [nonnegative]. Such matrices arise in various branches of mathematics and in many applications, including oscillations in mechanical systems~\cite{gk_book}, stochastic processes and approximation theory~\cite{karlin_tp}, planar resistor networks~\cite{curtis1998circular}, optimal allocation problems~\cite{bartroff2010optimal}, and many more~\cite{total_book,pinkus,fomin2000total}.

\updt{From here on we use~$n,m$ to denote integers that are larger than one.} One reason for the  importance of~TP matrices is their \emph{variation diminishing property}:
if~$A \in\R^{n\times m}$ is~TP, with~$n\geq m $,   then for  any~$x \in\R^m \setminus\{0\}$ the number of sign variations in~$Ax$ is smaller than or equal to the number of sign  variations in~$x$. 
Recently, it was shown that this property has important implications in the 
 asymptotic analysis of time-varying linear and 
nonlinear dynamical systems~\cite{margaliot2019revisiting,CTPDS,tpds_cdc,rola,eyal_k_posi,rola_ieee_tac}. In   these dynamical systems  the number of sign variations in the vector of derivatives
 is an integer-valued Lyapunov function.

Oscillatory matrices were introduced in the seminal work of Gantmacher and~Krein~\cite{gk_book}  who studied small vibrations of mechanical  systems. 
 %%%%
A matrix~$A\in\R^{n \times n}$ is called oscillatory if~$A$ is TN and there exists \updt{ a positive integer~$k$} such that~$A^k$ is TP. 
Thus, oscillatory matrices  are intermediary between TN and TP matrices.
Oscillatory matrices enjoy many special properties~\cite{gk_book}. For example, the eigenvalues of oscillatory matrices are real, positive and distinct, and the corresponding eigenvectors satisfy a special  sign pattern. 

An oscillatory matrix~$A$ must be  non-singular, as~$\det(A^k)>0$. 
Two useful characterizations of oscillatory matrices are the following.
\begin{Proposition}\label{prop:osc_GK}
\updt{\cite[Ch. 2]{total_book}} Let~$A\in\R^{n\times n}$  be a~TN matrix.
 Then~$A$ is oscillatory if and only if it is non-singular, 
 and~$a_{i,i+1}>0,a_{i+1,i}>0$, for all~$i=1,\dots,n-1$. 
\end{Proposition}

\begin{Proposition}\label{prop:osc_n-1}
\updt{\cite{classes}} A matrix~$A\in\R^{n \times n}$  is oscillatory if and only if~$A$ is TN, non-singular and irreducible.
\end{Proposition}

%We note in passing that zero entries and zero minors of oscillatory matrices are not arbitrary in nature; on the contrary, they ``throw a shadow" (see~\cite[Proposition 1.15]{pinkus} for more details).

  The exponent  of an oscillatory matrix~$A\in\R^{n\times n}$, denoted by~$r=r(A)$, is the \emph{least} positive integer such that~$A^{r}$ is TP. 
	It is well-known that if~$A$ is oscillatory then~$A^{n-1}$ is
	TP~\cite[Ch. 2]{total_book}, so~$r(A)\le n-1$. Yet, 
	deriving a closed-form expression for~$r$ for various classes of matrices
	is a nontrivial problem. 
%Note that since~$A$ is I-TN,~$A^v$ is TP, for all~$v\ge r$.

Oscillatory matrices have found  many applications~\cite{gk_book,price1968monotone,total_book,kardell2010total}. Recently, \updt{Katz et al.}~\cite{rami_osci} introduced the notion of \emph{oscillatory discrete-time systems} and used it in the analysis of certain discrete-time, time-varying nonlinear  systems. It was shown that if the mapping defining  the nonlinear dynamics is~$T$-periodic then any trajectory of the system either leaves any compact set or converges to a subharmonic trajectory, i.e. a trajectory that is periodic with a period of~$mT$. The \updt{positive integer~$m$}    is bounded by the exponent of an oscillatory matrix. 

\updt{Fallat and Liu}~\cite{fallat2007class} identified classes of oscillatory matrices with~$r(A)=n-1$. Motivated by the work in~\cite{fallat2007class,fallat2004remark}, we determine~$r(A)$ explicitly 
for several classes of oscillatory matrices
 (see \updt{Theorem~\ref{thm:deg} and Corollary~\ref{coro:main_explicit}}),
and   provide nontrivial upper bounds on~$r(A)$ 
for other classes (see Corollary~\ref{corr:bound}).

The remainder of this paper is organized as follows. The next section reviews known tools and results that will be used later on. Section~\ref{sec:main} describes our main results.
% Section~\ref{sec:proofs} includes the proofs of these results. 
Section~\ref{sec:conc} concludes and discusses possible directions for future research. 

We use standard notation. Vectors [matrices] are denoted by small [capital] letters. 
$\R^n$ [$\R_+^n$] is the set of vectors with~$n$ real [real and nonnegative] entries. 
%For any two integers $i$ and $j$ satisfying $i\le j$, $[i,j]:=\{i, i+1,\dots,j\}$.
For a matrix~$A\in\R^{n\times m}$, $a_{ij}$ or~$(i,j)$ denotes the entry of~$A$ in row~$i$ and column~$j$, and~$A^T$ is the transpose of~$A$.
%and $\tr(A)$ denotes the trace of~$A$. 
%Let~$\R^n_{++}:=\{v\in\R^n : v_i>0, i=1,\dots,n\}$, i.e. the set of all~$n$-dimensional vectors with positive entries. 
%We denote by~$A\gg0$ [$A\ge0$] a matrix with all entries positive [nonnegative].
 The square identity matrix is denoted by~$I$, with dimension that should be clear from the context.

\section{Preliminaries}
% %%%%%%%%%%%%%
We begin by reviewing notations and results that will be used later on. 
Let~$A\in\R^{n\times m}$.   Pick~$k_1 \in\{1,\dots,n\}$ and~$k_2\in\{1,\dots,m\}$, and let~$\alpha$ [$\beta$] denote a set of~$k_1$ [$k_2$] integers $1\leq i_1<\dots<i_{k_1}\leq n$
[$1\leq j_1<\dots<j_{k_2}\leq m$]. Then~$A[\alpha|\beta]$ denotes the~$k_1\times k_2$ sub-matrix of~$A$ containing the rows indexed by~$\alpha$ and the columns indexed by~$\beta$. When~$k_1=k_2$, 
we let~$A(\alpha|\beta):=\det(A[\alpha|\beta])$, that is, 
 the minor of~$A$ corresponding to the rows indexed by~$\alpha$ and columns indexed by~$\beta$.
 A minor corresponding to the same set of row and column indexes (i.e.~$A(\alpha|\alpha)$) is called a \emph{principal minor}. 

%Let~$A\in\R^{n\times m}$. Pick~$k \in\{1,\dots,\min\{n,m\}\}$, and let~$\alpha$ [$\beta$] denote a set of~$k$ integers $1\leq i_1<\dots<i_{k}\leq n$ [$1\leq j_1<\dots<j_{k}\leq m$]. Then~$A[\alpha|\beta]$ denotes the~$k\times k$ sub-matrix of~$A$ containing the rows indexed by~$\alpha$ and the columns indexed by~$\beta$. The determinant of~$A[\alpha|\beta]$, i.e. the minor of~$A$ corresponding to the rows indexed by~$\alpha$ and columns indexed by~$\beta$ is denoted by~$A(\alpha|\beta)$. A minor corresponding to the same set of row and column indexes (i.e. $A(\alpha|\alpha)$) is called a \emph{principal minor}. 
%For example, for~$A=\begin{bmatrix} 1&4&7\\3&5&6\\2&9&8\end{bmatrix}$, $\alpha=\{1,3\}$ and~$\beta=\{2,3\}$, we have $A(\alpha|\beta)=\det( \begin{bmatrix} 4& 7  \\ 9&8 \end{bmatrix})= -31$, $A(\alpha|\alpha)=\det( \begin{bmatrix} 1& 7  \\ 2&8 \end{bmatrix})= -6$, and $A(\beta|\beta)=\det( \begin{bmatrix} 5& 6  \\ 9&8 \end{bmatrix})= -14$. 

Pick~$A\in\R^{n\times m}$, $B\in \R^{m\times p}$,  
 and let~$C:=AB$.
The Cauchy-Binet formula~\cite[Ch.~1]{total_book} asserts that for any two
sets~$\alpha\subseteq\{1,\dots,n\}$, $\beta \subseteq \{ 1,\dots,p\}$, 
with the same cardinality~$k\in\{1,\dots,\min\{ n,m,p\}\}$, (i.e.~$|\alpha|=|\beta|=k$) we have 
\be\label{eq:cbfor}
				 C(\alpha|\beta)= \sum_{  |\gamma| =k}  A(\alpha|\gamma)  B(\gamma |\beta),
\ee
where the sum is over all~$\gamma=\{i_1,\dots,i_k\}$, with~$1\leq i_1<\dots<i_k\leq m$.  
Thus, every minor of~$AB$ is the sum of products of minors of~$A$ and~$B$.
%%
%For example, for~$n=m=p$ and~$k=n$, Eq.~\eqref{eq:cbfor} gives the well-known formula~$\det(AB)=\det(A)\det(B)$.
%%%
Note that~\eqref{eq:cbfor} implies that if~$A$ and~$B$ are both~TP [TN] then~$AB$ is~TP [TN], and that the product of an invertible~TN and~TP is~TP. In particular, this implies that if~$A\in\R^{n\times n}$ is oscillatory, then~$A^k$ is TP for all~$k\ge r(A)$.

Given~$A\in\R^{n\times n}$ and~$p\in\{1,\dots,n\}$, 
%Recall that each such  minor is defined by a set of row indexes~$1\leq i_1<i_2<\dots<i_p\leq n$ and column indexes~$1\leq j_1<j_2<\dots<j_p\leq n$, and is denoted by~$A(\alpha|\beta)$, where~$\alpha:=\{i_1,\dots,i_p\}$ and~$\beta:=\{j_1,\dots,j_p\}$.
the~$p$th \emph{multiplicative compound~(MC)} of~$A$, denoted~$A^{(p)}$,
 is the~$\binom{n}{p}\times  \binom{n}{p}$ matrix
that includes all the~$p\times p$ minors of~$A$ ordered lexicographically.
 The Cauchy-Binet formula yields
%\be\label{eq:miucpxdw}
$
(AB)^{(p)}=A^{(p)}B^{(p)},
$
%\ee
 justifying the term multiplicative compound. In particular,~$(A^k)^{(p)}=(A^{(p)})^k$  for all~$k\ge1$.

%Note that if~$A$ is TN and non-singular, then by Proposition~\ref{prop:TN_principal} all the diagonal elements of~$A^{(p)}$ are strictly positive, for all~$p\in\{1,\dots,n\}$. 

\updt{
An \emph{upper-right  corner minor}  of a matrix~$A\in\R^{n\times m}$ is a minor~$A(\alpha|\beta)$, where~$\alpha$ consists of the first~$k$ indexes and~$\beta$ consists of the last~$k$ indexes, for some $k\in\{1,\dots,\min\{n,m\}\}$, that is, 
the  upper-right corner minors are
\[
A(\{1,\dots,j\}|\{n-j+1,\dots,n\}) ,\quad   j=1,\dots,n.
\]
 A minor~$A(\alpha|\beta)$ is  a 
\emph{lower-left corner minor}
if~$A(\beta|\alpha)$ is an upper-right corner minor.
A~\emph{corner minor} is one that is either a lower-left or an upper-right corner minor.  
%%%%
 If~$A\in\R^{n\times n}$ then the corner minors of~$A$ are  the entries in the first row and last column, and the  last row and first column of every~$A^{(k)}$, 
i.e. $A^{(k)}_{1,\binom{n}{k}}$ and $A^{(k)}_{\binom{n}{k},1}$,  $k=1,\dots,n$, respectively}.

\begin{comment}	
\begin{table}[h]
\centering
\begin{tabular}{c c}
\hline\hline
Upper-right & Lower-left \\ [0.2ex] % inserts table %heading
\hline\hline
$A(\{1\}|\{n\})$ & $A(\{n\}|\{1\})$ \\
\hline
$A(\{1,2\}|\{n-1,n\})$ & $A(\{n-1, n\}|\{1, 2\})$ \\
\hline
$\vdots$ & $\vdots$   \\
\hline
$A(\{1,2,\dots,n-1\}|\{2,3,\dots,n\})$ & $A(\{2,3,\dots,n\}|\{1,2,\dots,n-1\})$ \\
\hline
\end{tabular}
\label{table:nonlin}
\caption{All upper-right and lower-left corner minors of~$A\in\R^{n\times n}$ (except for~$\det(A)$).}\label{tbl:corner}
\end{table}
\end{comment}

A matrix~$A\in\R^{n\times n}$ is~TP if  all the~$\sum_{p=1}^n \binom{n}{p}^2$ minors of~$A$ are positive. If~$A$ is known to be~TN then verifying  that 
a small subset of the minors are positive implies that~$A$ is~TP.  This is stated in
 the following result.  
%%%%%%(see, e.g.,\updt{~\cite[Theorem~4.3]{gasca1992total}} and~\cite[Theorem~3.1.10]{total_book}).
\begin{Proposition}\label{prop:corner}
\updt{(\cite[Ch. 3]{total_book},\cite{gasca1992total})}  Suppose that~$A\in\R^{n\times n}$
 is~TN. Then~$A$ is~TP if and only if
 all corner minors of~$A$ are positive.
\end{Proposition}

Our analysis of oscillatory matrices is motivated
 by the work of Fallat and Liu~\cite{fallat2007class}, 
and is based on the powerful successive elementary bidiagonal factorization of invertible TN matrices~\cite{Whitney1952,cryer_factorization},
 and their associated planar networks. 

%%%%%%%%%%%%%%%%%%%%%%%%%%%%%%%%%
\subsection{Successive Elementary Bidiagonal Factorization}
%%%%%%%%%%%%%%%%%%%%%%%%%%%%%%%%%
 Let~$E_{i,j}\in\R^{n \times n}_+$ denote  the matrix whose only nonzero entry is a one
in row~$i$ and column~$j$. For~$q\in\R$ and~$i\in\{2,\dots,n\}$,
 let~$L_i(q) := I + qE_{i,i-1}$ and~$U_i(q) := (L_i(q))^T$. For example, for~$n=3$, $L_2(4)=\begin{bmatrix} 1 & 0 & 0 \\ 4 & 1 & 0 \\ 0 & 0 & 1 \end{bmatrix}$. The matrices~$L_i(q)$ and~$U_i(q)$ are called \emph{elementary bidiagonal~(EB) matrices}.
Several useful relations of EB matrices are:
\begin{alignat}{2} \label{eq:Lij_flip}
%%%%%%%%%%%%%%%%%%%%%%%%%
L_i(0)&=U_i(0)=I,\quad \;\; && i=2,\dots,n,  \nonumber \\
L_i(x)^{-1}&=L_i(-x),  &&i=2,\dots,n,  \nonumber \\
U_i(x)^{-1} &= U_i(-x),  && i=2,\dots,n, \nonumber \\
L_i(x) L_j(y) &= L_j(y) L_i(x) ,  &&  |i-j|>1, \nonumber \\
U_i(x) U_j(y) &= U_j(y) U_i(x),   &&|i-j|>1, \nonumber \\
L_i(x) U_j(y) &= U_j(y) L_i(x) ,  && i \ne j.
\end{alignat}

 \updt{
\begin{Theorem}\label{prop:SEB}
\cite[Ch. 2]{total_book} Let~$A\in\R^{n\times n}$  be  an invertible TN~(I-TN) matrix. Then~$A$ can be factorized as:
\begin{align}\label{eq:EB}
A =& [L_n(\ell_1)\cdots L_2(\ell_{n-1})][L_n(\ell_{n})\cdots L_3(\ell_{2n-3})]\cdots [L_n(\ell_k)]D \nonumber \\
& [U_n(u_k)][U_{n-1}(u_{k-1})U_n(u_{k-2})]\cdots[U_2(u_{n-1})\cdots U_n(u_1)],
\end{align}
where~$k:=(n-1)n/2$,
 $\ell_i, u_i\ge0$, $i=1,\dots,k$, and $D\in\R^{n\times n}$ is a diagonal matrix with positive   entries.
\end{Theorem}
}

\begin{Remark}
Note that~$k\geq n-1$ for all~$n\geq 2$. 
The representation~\eqref{eq:EB}  is called the~{successive elementary bidiagonal~(SEB) factorization} of the I-TN matrix~$A$, and this factorization is unique~\cite[p.~53]{total_book}.  The~$\ell_j$'s and~$u_k$'s are called the {multipliers} in the factorization. 
\end{Remark}

For example, for~$n=4$ we have~$k=6$, 
so any I-TN matrix~$A\in\R^{4\times4}$ can be factorized as:
\begin{align*}
A=&[L_4(\ell_1)L_3(\ell_2)L_2(\ell_3)][L_4(\ell_4)L_3(\ell_5)][L_4(\ell_6)]
\\&D[U_4(u_6)][U_3(u_5)U_4(u_4)][U_2(u_3)U_3(u_2)U_4(u_1)].
\end{align*}

The derivation of~\eqref{eq:EB} is based on the well-known Neville elimination 
process~\cite{gasca1992total}.  
 The following example demonstrates this.  
We use~$\diag(d_1,\dots,d_n)$ to denote the~$n\times n$ 
diagonal matrix with diagonal entries~$d_1,d_2,\dots,d_n$.
%%%
\begin{Example}
Consider the matrix
\be
A=\begin{bmatrix} 1 & 1 & 1 \\ 
                              2 & 4 & 8 \\
                              2 & 10 & 29
                              \end{bmatrix}.
\ee
It is straightforward to verify that~$A$ is~I-TN.                              
The first step in the Neville elimination process is to use the second row to null entry~$(3,1)$ in~$A$. This is done by multiplying~$A$  from the left  by~$L_3(-1)$  yielding
\[
L_3(-1)A =\begin{bmatrix}
1 & 1 & 1 \\
2 & 4 & 8 \\
0 & 6 & 21
\end{bmatrix}.
\]
We next use the first row to null entry~$(2,1)$ in~$L_3(-1)A$ using~$L_2(-2)$:
\[
L_2(-2)L_3(-1)A=\begin{bmatrix}
1 & 1 & 1 \\
0 & 2 & 6 \\
0 & 6 & 21
\end{bmatrix}.
\]
Next, we use the second row to null entry~$(3,2)$ in~$L_2(-2)L_3(-1)A$ using~$L_3(-3)$:
\[
L_3(-3)L_2(-2)L_3(-1)A=\begin{bmatrix}
1 & 1 & 1 \\
0 & 2 & 6 \\
0 & 0 & 3
\end{bmatrix} :=P.
\]
Applying similar row operations to~$P^T$ yields~$L_3(-2)L_2(-1)L_3(-1)P^T=D$, 
where~$D:=\diag(1,2,3)$. Thus,~$L_3(-2)L_2(-1)L_3(-1)A^TL_3(-1)^TL_2(-2)^TL_3(-3)^T=D$, and using~\eqref{eq:Lij_flip} yields
\be\label{eq:exp1}
A = [L_3(1)L_2(2)][L_3(3)]D[U_3(2)][U_2(1)U_3(1)],
\ee
which is the SEB factorization of~$A$.
\end{Example}

TP and oscillatory matrices can be characterized in terms of their SEB factorization.
%%%%
\begin{Proposition}\label{prop:seb_TP}
\updt{\cite[Ch. 2]{total_book}}  A matrix~$A\in\R^{n \times n}$  is TP if and only if in the SEB factorization~\eqref{eq:EB}~\updt{all the multipliers are positive}.
\end{Proposition}

\begin{Proposition}\label{prop:seb_OSC}
\updt{\cite[Ch. 2]{total_book}}  A matrix~$A\in\R^{n \times n}$ is oscillatory if and only if in the SEB factorization~\eqref{eq:EB} at least one of the multipliers from each of~$L_n,L_{n-1},\dots,L_2$, and from each of~$U_n,U_{n-1},\dots,U_2$ is positive.
\end{Proposition}

The case where exactly one  of the multipliers above is positive defines a \emph{basic oscillatory} matrix. 
\begin{Definition}
\updt{\cite[Ch. 2]{total_book}}  An I-TN matrix~$A\in\R^{n \times n}$  is called  \emph{basic oscillatory} if and only if in the SEB factorization~\eqref{eq:EB} {exactly} one of the multipliers from each of~$L_n,\dots, L_2$, and exactly one from each  of~$U_n,\dots,U_2$ is positive.
\end{Definition}
%%%%%%%%%%%%%%%
For example, the factorization of 
\[
A = 
\begin{bmatrix}
1 & 6 & 0 & 0 \\
2 & 13 & 4 & 20 \\
2 & 13 & 5 & 25 \\
0 & 0 & 3 & 16
\end{bmatrix},
\]
is~$A=L_3(1)L_2(2)L_4(3)IU_3(4)U_4(5)U_2(6)$, so~$A$ 
is a basic oscillatory matrix. Basic oscillatory matrices may be viewed as minimal oscillatory matrices in the sense of the number of SEB factors involved. 
 Every oscillatory matrix~$A$ can be written in the form $A =
SBT$, where~$B$ is a basic oscillatory matrix and $S, T$ are  I-TN~\cite[Ch.~2]{total_book}.
%%
%%
%The following result shows that any oscillatory matrix is a product of a basic oscillatory and I-TN matrices.
%\begin{Proposition}
%\updt{\cite[Ch. 2]{total_book}}  Any oscillatory matrix~$A\in\R^{n\times n}$  
% can be written in the form~$A=A_1BA_2$, where~$B\in\R^{n\times n}$ is basic oscillatory and both~$A_1,A_2\in\R^{n\times n}$  are I-TN.
%\end{Proposition}

Let~$A\in\R^{n\times n}$  be oscillatory. 
\updt{Fallat and Liu}~\cite{fallat2007class} derived a necessary and sufficient condition for its exponent to be~$n-1$ using  the SEB factorization. 
For example, let
\be\label{eq:jacobi}
A=L_2(\ell_1)L_3(\ell_2)\cdots L_n(\ell_{n-1})DU_n(u_{n-1})\cdots U_3(u_2)U_2(u_1),
\ee
where~\updt{every~$\ell_i$ and $u_i$  is  positive}. Then~$A$ is a \emph{tridiagonal} basic oscillatory matrix (also referred to as a Jacobi matrix), and it is not difficult to see that the entries~$(n,1)$ and~$(1,n)$ in~$A^{n-2}$ are zero. Thus,~$r(A)>n-2$ and since~$r(A)\le n-1$, $r(A)=n-1$. In addition, the cases where
\be\label{eq:bosc1}
A=L_n(\ell_1)L_{n-1}(\ell_2)\cdots L_2(\ell_{n-1})DU_2(u_{n-1})\cdots U_{n-1}(u_2)U_n(u_1),
\ee
where~\updt{every~$\ell_i$ and $u_i$ is   positive}, also yield~$r(A)=n-1$ \updt{(see~\cite{fallat2007class})}.

More generally, \updt{Fallat and Liu}~\cite{fallat2007class} established conditions on~$\ell_i, u_i$, $i=1,\dots,k$, in~\eqref{eq:EB} that yield~$r(A)=n-1$, and conditions that yield~$r(A)\le n-2$. 
For example, \updt{for~$n=4$, pick~$L = L_4(\ell_1) L_2(\ell_2) L_3(\ell_3)$ or~$L = L_3(\ell_1) L_2(\ell_2) L_4(\ell_3)$, and~$U = U_4(u_1) U_2(u_2) U_3(u_3)$ or~$U = U_3(u_1) U_2(u_2) U_4(u_3)$}, 
%For example, for~$n=4$, pick one form of~$L$ from:
%\begin{align*}
%L &= L_4(\ell_1) L_2(\ell_2) L_3(\ell_3), \\
%L &= L_3(\ell_1) L_2(\ell_2) L_4(\ell_3),
%\end{align*} 
%and one form of~$U$ from
%\begin{align*}
%U &= U_4(u_1) U_2(u_2) U_3(u_3), \\
%U &= U_3(u_1) U_2(u_2) U_4(u_3),
%\end{align*} 
with all \updt{$\ell_i,u_i$ positive}. Then~$A=LDU$,
where~$D$ is diagonal with positive diagonal entries,
 implies that~$A^2$ is TP, so~$r(A)=2$ (note that~$A$ is not TP). On the other hand, $A=L_2(\ell_1)L_3(\ell_2)L_4(\ell_3)DU$ or~$A=L_4(\ell_1)L_3(\ell_2)L_2(\ell_3)DU$, with~\updt{$\ell_i$ positive} and any~$U$, implies
 that~$A^2$   is not TP, so~$r(A)>2$ and thus~$r(A)=3$.

%%%%%%%%%%%%%%%%%%%%
\subsection{Planar Networks}
%%%%%%%%%%%%%%%%%%%%
An \emph{elementary weighted diagram} is a weighted and directed
graph consisting of~$n$ source vertices (on the left of the graph)
and~$n$ sink vertices (on the right). The  sources and
sinks    are numbered consecutively \emph{from bottom to top}. 
All edges in the diagram are directed from left to right.  
 
 \begin{figure*}
\begin{center}
\includegraphics[width= 11cm,height=5cm]{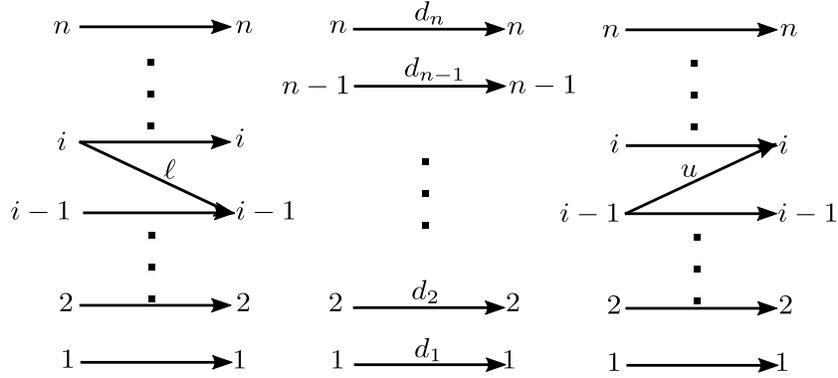}
\caption{Elementary diagrams associated with the EB matrix~$L_i(\ell)$ (left), 
\updt{diagonal matrix~$D$ (middle), and EB matrix~$U_i(u)$ (right)}. Unmarked edges have  weight   one.}\label{fig:EB_mats}
\end{center}
\end{figure*}

The elementary weighted diagrams of the EB matrices~$L_i(\ell)$, $U_i(u)$, and a diagonal matrix~$D=\diag(d_1,\dots,d_n)$ are depicted in Figure~\ref{fig:EB_mats}. All edges in the figure  are directed from left to right, and unmarked edges have weight   one. The source [sink] vertices on the left [right] side represent the matrix rows
 [columns] indexes, and an edge of weight~$q$ from row index~$i$ to column index~$j$ implies that 
the entry~$(i,j)$ in the associated matrix is equal to~$q$.

Given an elementary diagram,  pick~$k\in\{1,\dots,n\}$,~$1\leq i_1<\dots<i_k \leq n$ and~$1\leq j_1<\dots<j_k \leq n$. We define a \emph{family} of paths that connect the sources~$\{i_1,\dots,i_k\}$
to the sinks~$\{j_1,\dots,j_k\}$ as follows. 
Each family contains~$k$   \emph{vertex-disjoint} paths (i.e. paths that are non-intersecting and non-touching) joining the vertices~$\{i_1,\dots,i_k\}$ on the left side of the diagram with the vertices~$\{j_1,\dots,j_k\}$ on the right side. 
The \emph{weight of a path} is defined to be the product of the weights of its edges. Note that in an elementary diagram a path contains a single edge, but in the more general diagrams defined below a path   consists of several edges. 
 The \emph{weight of the family} is defined to be the product of the weights of its~$k$ paths. 

%%%%%%%%%%%
\begin{Example}\label{exa:loep}
%%%%%%%%%%%%%%%%%
Consider the diagram on the left-hand side of Figure~\ref{fig:EB_mats}, 
corresponding to~$Q=L_i(\ell)$, 
and consider the sources~$\{ i-1 , i \}$ and the sinks~$\{i-1 , i\}$.  
There is one corresponding family, namely,~$\{ i-1 \to i-1, i \to i\}$ (the path~$i\to i-1$ cannot be included, as the   paths must be vertex-disjoint). 
The weight  of each of the two paths is one, and so the weight of the family  
 is also one. 

As another example, suppose that~$Q=L_i(\ell)$,  with~$i>2$,
and consider the sources~$\{ 1 , i \}$ and the sinks~$\{1 , i-1\}$.  
There is one corresponding family, namely,~$\{ 1 \to  1, i \to i-1\}$.
The weight  of   the first [second] path is~$1$ [$\ell$], and so the weight of the family  
 is~$\ell$. 
%%%%%%%%%%%%%%%%%%%
\end{Example}

We now review  an   
  important property  of these diagrams~\cite{total_book}. 
 Let~$Q\in\R^{n\times n}$  be a matrix represented by any one of the diagrams in~Figure~\ref{fig:EB_mats}. 
For any~$k\in\{1,\dots,n\}$ and any~$1\leq i_1<\dots<i_k \leq n$ and~$1\leq j_1<\dots<j_k \leq n$
consider all the families of~$k$ vertex-disjoint paths 
 joining the vertices~$\{i_1,\dots,i_k\}$ on the left side of the diagram with the vertices~$\{j_1,\dots,j_k\}$ on the right side. This set  of families is unique. 
We define the \emph{weight of the set of families} as the sum
of the weights of each family in the set. 
 Then the minor
\[
Q(\{i_1,\dots,i_k\}|\{j_1,\dots,j_k\})
\]
is equal to weight of the set of families (and is zero if and only if 
there is no such family). 
For example, consider again 
the diagram on the left-hand side of Figure~\ref{fig:EB_mats}, corresponding to~$Q=L_i(\ell)$.
Then
\updt{ % here I changed \begin{align} to \begin{align*}
%\begin{align*}
\[
Q(\{i-1,i \} | \{i-1,i\}) =\det( \begin{bmatrix}  1&0\\\ell & 1 \end{bmatrix} )=1,
%\end{align*}
\]
and
%\begin{align*}
\[
Q(\{ 1,i \} | \{1,i-1\}) =\det( \begin{bmatrix} 
 1&0\\
0&\ell   \end{bmatrix} )  \\
=\ell,
%\end{align*}
\]
}
 and this agrees with the results
 in Example~\ref{exa:loep}.

\begin{figure*}
\begin{center}
\includegraphics[width= 12cm,height=5cm]{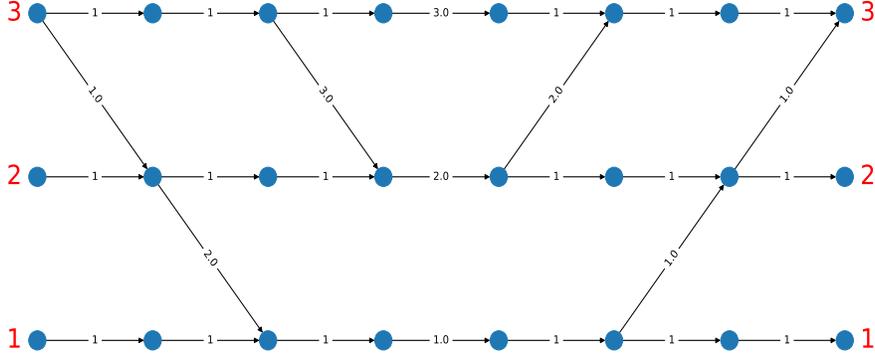}
\caption{Planar network associated with the I-TN matrix~$A$ in~\eqref{eq:exp1}.}\label{fig:EB_exp1}
\end{center}
\end{figure*}

Given the factorization~\eqref{eq:EB}, the corresponding~\emph{planar network} associated with~$A$ is obtained by concatenating (in order) left to right the elementary diagrams associated with the EB matrices and~$D$ in~\eqref{eq:EB}. For example, Figure~\ref{fig:EB_exp1} depicts the planar network associated with the matrix~$A$ in~\eqref{eq:exp1}.
 
We note in passing that, conversely, any planar network can be associated with a TN matrix~\cite{lindstrom1973vector}. This association is unique in case of an I-TN matrix. 

\begin{Remark}\label{remk:planar}
Since the planar network of any I-TN matrix is associated with an SEB factorization, it admits a special structure.  First, all horizontal edges have positive weights, where horizontal edges corresponding to the matrix~$D$ (in the center of the network) have weight~$d_{i}$, $i=1,\dots,n$, and all other horizontal edges have weight one. In addition, the left-hand side [right-hand side] of the network, corresponding to the product of all~$L_i$ [$U_i$] matrices, consists of only downward [upward] pointing  diagonal edges. This structure holds for all I-TN matrices. 
Only the weights of the diagonal edges and the horizontal edges corresponding to the matrix~$D$ differ among different I-TN matrices.
%%%%%%%%%%%%%%%%%%%%%%%%%%%%%%%%%%%%%%%%%%%%%%%
\end{Remark}

  The following important result   associates the minors of~$A$  to its   planar network.
Its proof follows from 
 the Cauchy-Binet formula.

\begin{Proposition}\label{prop:minor_planar}
%%%%%%%%%%%%%%%%%%%%%%%%%%%%%%%%%%%%%%%%%%%%%
\updt{\cite[Ch. 2]{total_book}}  Let~$A:=A_1\cdots A_p$, where each~$A_i$ is either an EB matrix ($L$ or~$U$) or a diagonal matrix  with positive diagonal entries. Recall that the network associated with~$A$ is obtained by concatenating left to right the diagrams associated with~$A_1,\dots, A_p$, respectively. Then~$A(\{i_1,\dots,i_k\}|\{j_1,\dots,j_k\})$ is equal to the  weight of the   set of vertex-disjoint families in the network connecting the sources~$\{i_1,\dots,i_k\}$
	to the sinks~$\{j_1,\dots,j_k\}$. 
	In particular,~$a_{i j}=A(\{i\}|\{j\})$ is equal to the sum  of all the   weights of paths that 
	join  source~$i$  
	to   sink~$j$.
%%%%%
\end{Proposition}

\begin{figure*}
\begin{center}
\includegraphics[width=13cm,height=5cm]{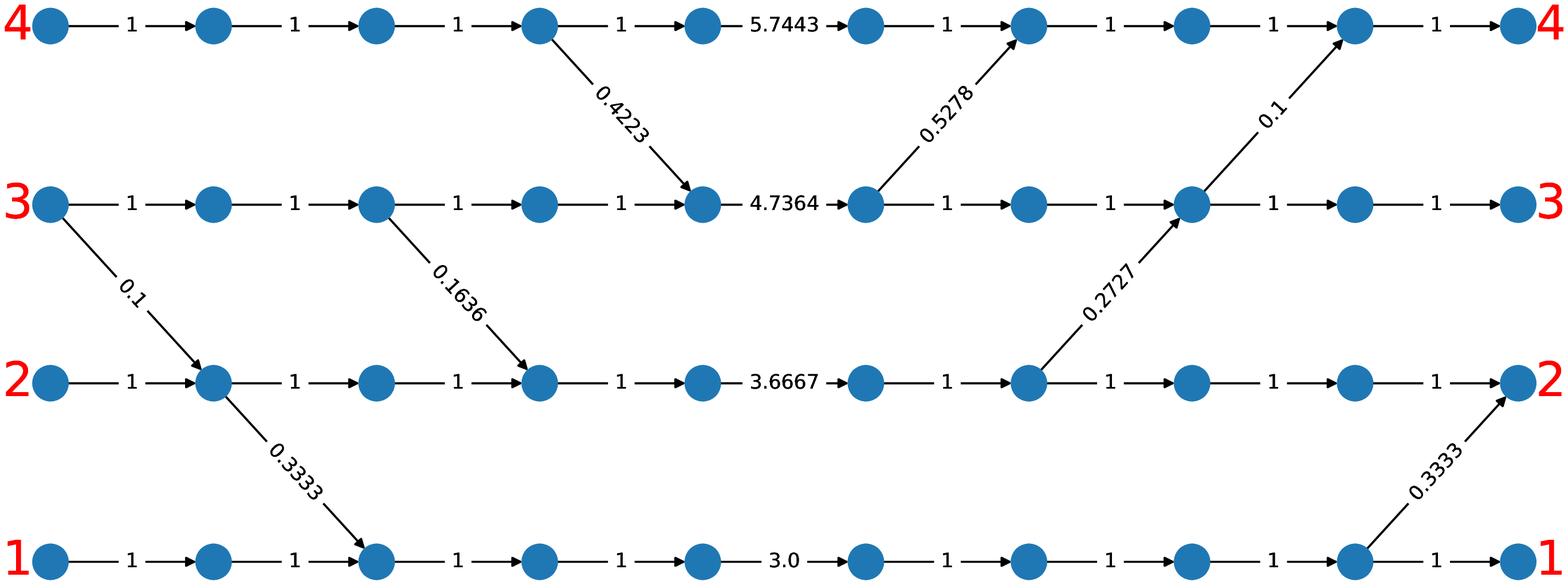}
\caption{Planar network associated with~$A$ in~\eqref{eq:A1}.}\label{fig:A1}
\end{center}
\end{figure*}

The next example demonstrates Proposition~\ref{prop:minor_planar}.
%%%%%%%%%%%%%%%%%%%%%%%
\begin{Example}
Let
\be\label{eq:A1}
A=\begin{bmatrix} 3 & 1 & 0 & 0 \\ 1 & 4 & 1 & 0.1 \\ 0.1 &1 & 5 & 3 \\ 0 & 0 & 2 &7 \end{bmatrix}.
\ee
It is straightforward to verify that~$A$ is I-TN. The associated network of~$A$ is depicted in Figure~\ref{fig:A1} (all numerical values  in this paper are to four-digit accuracy).

There is only one directed path joining the  source~$2$   with the  sink~$4$, and the weight of this path is~$3.6667\cdot0.2727\cdot 0.1\approx0.1$. Indeed, $a_{2,4}=0.1$. 

There are  two paths  that join source~$2$  with sink~$2$: the horizontal path with weight~$3.6667$, and the path with weight~$0.3333\cdot 3 \cdot 0.3333\approx 0.3333$. The sum of these two path weights is~$4$, and this is equal to~$a_{2 2}$. 

 As another example, consider~$A(\{1,3\}|\{2,3\})$. In this case, there are three families of vertex-disjoint paths joining the sources~$\{1,3\}$ 
 with the sinks~$\{2,3\}$:
\begin{enumerate}
\item The path joining source~$1$
 with sink~$2$ 
with weight~$3\cdot0.3333\approx 1$, and the path joining source~$3$ with sink~$3$  with weight~$4.7364$. This family weight is then~$1\cdot4.7364=4.7364$.
\item The path joining source~$1$ 
with sink~$2$  
 with weight~$3\cdot0.3333\approx1$, and the path joining source~$3$ with sink~$3$   with weight~$0.1\cdot3.6667\cdot0.2727\approx 0.1$. This family weight is then~$1\cdot0.1=0.1$.
\item The path joining source~$1$  
 with sink~$2$ with weight~$3\cdot0.3333\approx1$, and the path joining source~$3$ 
 with sink~$3$  with weight~$0.1636\cdot3.6667\cdot0.2727\approx0.1636$. 
This family weight is then~$1\cdot0.1636=0.1636$.
\end{enumerate}
The sum over these three families is~$4.7364+0.1+0.1636 =5$,
and~$A(\{1,3\}|\{2,3\})=\det(\begin{bmatrix} 1 & 0 \\ 1 & 5 \end{bmatrix})=5$.
\end{Example}

Proposition~\ref{prop:minor_planar} links the minors of an I-TN matrix  to the topology of the associated planar network. This has many theoretical and practical implications. For example, the following known proposition follows immediately from Proposition~\ref{prop:minor_planar}, as horizontal edges with positive weights always exist in \updt{the planar network of any I-TN matrix} (see Remark~\ref{remk:planar}). 
%As Proposition~\ref{prop:minor_planar} is a very important result as it links the minors of I-TN matrices to the (very regular) topology of the associated planar network. Indeed, it is one of the main tools we use to obtain our results below. One can appreciate the importance of this result, for example, by noting that the following known proposition follows immediately from Proposition~\ref{prop:minor_planar}, as horizontal edges with positive weights always exist in any planar diagram of I-TN matrices.

\begin{Proposition}\label{prop:TN_principal}
\updt{\cite[Ch. 1]{pinkus}} Let~$A\in\R^{n\times n}$
 be an I-TN matrix. Then every principal minor of~$A$ is  positive.
\end{Proposition}

Proposition~\ref{prop:minor_planar}  is one of the 
main tools we \updt{use} to derive the  results below.

%%%%%%%%%%%%%%%%%%%%%%%
\section{Main Results}\label{sec:main}
%%%%%%%%%%%%%%%%%%%%%%%%

\updt{ % =======================================================
Since an oscillatory matrix~$A$ is in particular~TN, Proposition~\ref{prop:corner}
implies that verifying that~$A^k$ is TP is equivalent to verifying that all
the corner minors of~$A^k$ are positive.
The next  result shows that when~$A$ is 
the product of matrices that appear in the SEB
factorization  it is possible  in some sense to ``decouple''
the examination of lower-left corner and upper-right  corner minors.
%%%%%%
 We say that a diagonal matrix~$D\in\R^{n\times n }$ is  
  \emph{positive} if all its diagonal entries are positive. 
%%%%%%%%%%%%%%%%%%%%%%%%%%%%%%%%%%
\begin{Proposition}\label{prop:LUcorner}
%%%%%%%%%%%%%%%%%%%%%%%%%%%%%%%%%%%%%%%%%%%%%%%%%%%%%%%%%%
Let~$A:=A_1\cdots A_p$, where each~$A_j\in\R^{n\times n}$
is either an EB matrix ($L_i$ or~$U_i$) or a positive diagonal matrix. 
Pick~$c\in\{1,\dots,n-1\}$, and let~$\alpha:=\{n-c+1,\dots,n\}$,   $\beta:=\{1,\dots,c\}$, 
so that~$A(\alpha|\beta)$ is a lower-left corner minor of~$A$.
 %%%%
Let~$w=w(A,c )$
be the minimal integer such that
in the planar network associated with~$A^w$ there is 
a path connecting~$ x$ to~$x-n+c$ for any~$x\in \alpha$. 
Then~$w$ does not change if we replace every~$U_i$ in~$A$ with any positive diagonal matrix.
%%%%%%%%%%%%%%% 
\end{Proposition}

In other words, the non-horizontal edges in the planar network 
corresponding to the~$U_i$s  do not affect~$w$. 

 By transposition, this implies a similar result for any
 upper-right corner minor~$A(\beta|\alpha)$,
namely, 
  the minimum number of copies~$v$ of~$A$ 
	needed to guarantee the existence of a path connecting~$y$ to~$y+n-c$
	for any~$y\in\beta$ does not depend on the non-horizontal edges in the~$L_i$s.

  We now describe another implication of
  Proposition~\ref{prop:LUcorner}.	Let~$A=L_aD_aU_a$,~$B=L_bD_bU_b$ 
	be I-TN 
	such that~$L_a$ and~$L_b$ [$U_a$ and~$U_b$] denote the product of all the~$L_i$
	[$U_i$] matrices in the SEB factorization of~$A$ and~$B$. 
Suppose that a parameter in the bidiagonal factorization of~$L_a$ is positive if and only if~(iff) the
corresponding parameter in the bidiagonal factorization of $L_b$ is positive.
 Then  
	a lower-left   corner minor of~$A$ is positive iff
	the same minor in~$B$  is positive.

%%%%%%%%%%%%%%%%%%%%%%%%%%%%%%%%%%%%%%%%%
{\sl Proof of Proposition~\ref{prop:LUcorner}.} 
%%%%%%%%%%%%%%%%%%%%%%%%%%%%%%%%%%%%%%%%%
If~$w(A,c)$ is infinite then clearly it will remain infinite if we replace 
 every~$U_i$ in~$A$ with some positive diagonal matrix.
 So we may assume that~$w(A,c)$ is finite.

Let~$G^w$ denote 
the network corresponding to~$w$ concatenated copies of the network of~$A$.
Then~$w(A,c)>0$ iff~$G^w$ includes 
  a family of
\emph{vertex disjoint} paths that include: 
a path~$P_c$ from~$n$ to~$c$;
a path~$P_{c-1}$ from~$n-1$ to~$c-1$; $\dots$; and
a path~$P_1$ from~$n-c+1$ to~$1$. 
If none of these~$P_i$s includes
an   upward pointing diagonal edge from some~$U_s$ then we are done.
Otherwise, there exists a \emph{minimal}~$i$ such that~$P_i$
includes an   upward pointing diagonal edge from some~$U_s$.
This edge points from a node~$v_r $ to a node~$v_{r+1}$ for some~$r\in\{1,\dots,n-1\}$. 
The continuation of~$P_i$ after this edge 
must include  a downward pointing diagonal edge from a node~$v_{r+1}$ to~$v_r $.
Let~$\tilde P_i$ be the path obtained from~$P_i$
by replacing this part in~$P_i$ by a set of
   horizontal arcs connecting~$v_r$ to~$v_r$. Then~$\{ P_c,P_{c-1},\dots,P_{i+1},\tilde P_i, P_{i-1}, \dots, P_1\}$ is a set of vertex disjoint paths that connect every~$x\in \alpha$ to~$ (x-n+c)\in \beta$ using  the same number of copies of~$A$.  Thus, we may assume that~$P_i$ does not include any 
	upward pointing diagonal edges. Now we can apply the same 
	idea to any upward pointing diagonal edge in~$P_{i+1}$ and so on.~\hfill{$\square$}
} % end of \updt

To proceed, 
we introduce a notation  that simplifies
 the representation of the factorization
 in~\eqref{eq:EB}. 
 For~$i\in\{2,\dots,n\}$, 
define:
\[
x_i:=1+(i-2)n-\sum_{j=1}^{i-2}j= (i-2)n-\frac{(i-3) i}{2},
\]
and
\begin{align}\label{eq:l^i}
\ell^i  :=\begin{bmatrix} \ell_{x_i} & \cdots & \ell_{x_i+n-i} \end{bmatrix}^T, \text{ and }
u^i  :=\begin{bmatrix} u_{x_i} & \cdots & u_{x_i+n-i} \end{bmatrix}^T.
\end{align}
Thus,
\begin{align*}
\ell^2&=\begin{bmatrix} \ell_1 & \cdots & \ell_{n-1} \end{bmatrix}^T\in\R_+^{n-1}, \\
\ell^3&=\begin{bmatrix} \ell_n & \cdots & \ell_{2n-3} \end{bmatrix}^T\in\R_+^{n-2},\\
&\vdots \\
\ell^n &= \ell_k\in\R_+,
\end{align*}
and similarly for~$u^i$.
Let~$\ell^i_j$ [$u^i_j$] 
 denote  the~$j$'th entry in~$\ell^i$ [$u^i$], and
let
\begin{align}\label{eq:WQ}
W_i(\ell^i)  &:= \prod_{j=n}^{i} L_j(\ell^i_{n-j+1}), \nonumber \\
Q_i(u^i) &:= \prod_{j=i}^n U_j(u^i_{j-i+1}). 
\end{align}
Then the SEB factorization~\eqref{eq:EB} can be written more succinctly as
\be\label{eq:EB_WQ}
A = W_2(\ell^2)\cdots W_n(\ell^n) D Q_n(u^n)\cdots Q_2(u^2).
\ee

\updt{
%%%%%%%%%%%%%%%%%%

For~$A\in \R^{n\times n }$, let~$r_\ell(A)$ denote the minimal integer~$w>0$
such that all the lower-left corner minors of~$A^{w}$, except perhaps for~$\det(A)$, are positive. Similarly, let~$r_u(A)$ denote the minimal integer~$v>0$
such that all the upper-right corner minors of~$A^{v}$, except perhaps for~$\det(A)$, are positive.
%%%%%%%%%%%%%%%%%%%%%%%%%%%%%%%%%%%
\begin{Theorem}\label{thm:deg}
%%%%%%%%%%%%%%%%%%%%%%%%%%%%%%%%%%%%%%
Let~$A\in \R^{n\times n}$ be oscillatory, and write 
its SEB factorization as in~\eqref{eq:EB_WQ}.
Then
\[
r(A)=\max\left \{r_\ell \left (W_2(\ell^2)\dots W_n(\ell^n)\right ),r_u\left (Q_n(u^n)  \dots Q_2(u^2)\right ) \right \}.
\]
%%%%%%%%%
\end{Theorem}
%%%%%%%%%%%%%%%%%%%
{\sl Proof of Theorem~\ref{thm:deg}.}
%%%%%%%%%%%%%%%%%%%%%%%%%%%%%%%%%
Since~$A$ is oscillatory, it is~TN and~$\det(A)>0$. Hence,~$r(A)=\max\left \{r_\ell \left (A\right ),r_u\left (A\right ) \right \}$. 
Combining this with the fact that every~$W_i$  [$Q_i$] is a product of~$L_j$s
[$U_j$s] and Proposition~\ref{prop:LUcorner} completes  the proof.~\hfill{$\square$} 

Theorem~\ref{thm:deg} can be applied to determine~$r(A)$ explicitly for 
large classes of
oscillatory matrices. Indeed, suppose that for some specific structure of~$W_2(\ell^2)\cdots W_n(\ell^n)$, denoted~$W$,   we can explicitly  
determine~$r_\ell\left (W\right )$. A typical case is when there exists a ``worst-case'' lower-left corner minor~$(\alpha|\beta)$, that is, if~$W^k(\alpha|\beta)>0$ for some~$k>0$
 then \emph{all}   lower-left corner minors of~$W^k$ (except perhaps for~$\det(W^k)$) are positive. If it is possible to determine the minimal~$j$ such that~$W^j(\alpha|\beta)>0$ then~$r_\ell\left (W\right )=j$.

By transposition, we also know~$ r_u(W^T)$. 
Thus, we know~$r(A)$ for~$A=WDW^T$, for any positive diagonal matrix~$D$. 
To demonstrate these ideas, we 
 now define two classes of~$n\times n$ matrices, denoted~$Z_i(s)$, $i=1,2$.
Each class includes matrices parametrized by a non-negative integer~$s$,
 taking values in a domain~$D_i$, $i=1,2$. In   the definitions of these classes below
\emph{it is always assumed that all the multipliers are positive}.
%%%%%%%%%%%%%%%%%%%%%%%%%%%%%%%%%%%%%%%%%%%%%
\begin{Definition}
%%%%
$Z_1(s)$, $s\in D_1:=\{2,\dots,n\} $, is the class of matrices in the form:
\[
W_2(\ell^2)\cdots W_s(\ell^s) ,
\]
where all   the multipliers are positive.
%%%%%
%%%
\end{Definition}
The   planar network of~$ Z_1(s)$ is
 the product of the networks
 for~$W_2(\ell^2)$, $\dots$, $ W_{s}(\ell^{s})$. This includes horizontal lines,
 and~$s-1$ downward 
 pointing  diagonal lines: the first   
from~$n$ to~$1$, the second from~$n$ to~$2$, and so on.  Figure~\ref{fig:EB_Z1} depicts this 
network for~$n=6$ and~$s\in\{2,\dots,6\}$. 
%% (see Table~\ref{tbl:corner}). 
%It is clear that only the diagonal lines with an ``upper-left to lower-right orientation'' in the left-hand side of the associated planar network (and all the horizontal lines) can contribute to the family of vertex-disjoint paths.

\begin{figure*}
\begin{center}
\includegraphics[width= 13cm,height=7cm]{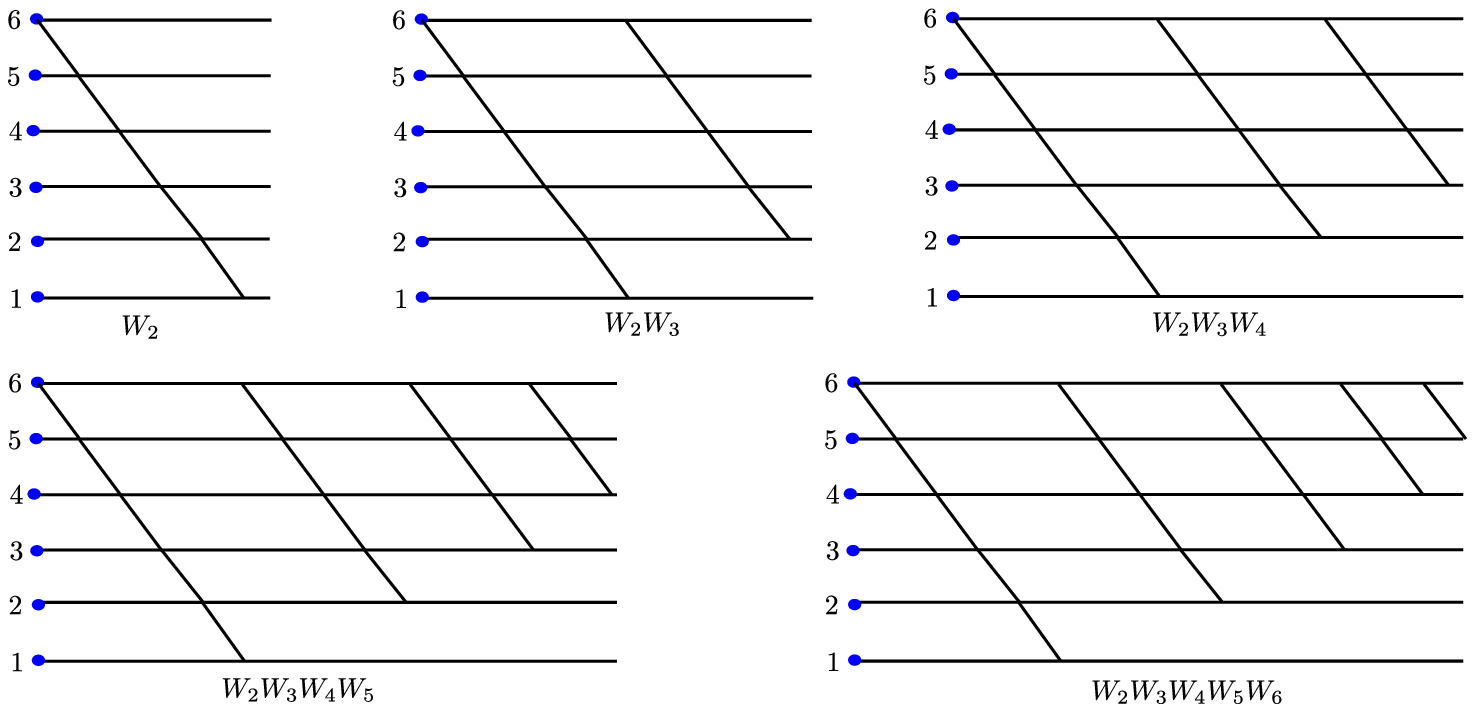}
\caption{\updt{Planar network  corresponding to~$Z_1(s)$  for~$n=6$ 
and~$s\in\{2,\dots,6\}$. 
Upper-left figure corresponds to~$s=2$, upper-middle to~$s=3$, 
upper-right to~$s=4$, lower-left to~$s=5$, and lower-right to~$s=6$.}}\label{fig:EB_Z1}
\end{center}
\end{figure*}

For~$y\in\R$, let~$\lceil y \rceil$ denote  the smallest  integer that is larger than or equal to~$y$.
%%%%%
\begin{Proposition}\label{prop:zrwq}
%%%%%%%%%%%%%%%%%%%%%
Pick~$s\in D_1$ and a  matrix~$A\in Z_1(s)$.  Let~$\mu_1:=\left\lceil  \frac{n-1}{s-1} \right\rceil$.
 Then~$r_\ell(A)=\mu_1 $. 
%%%%%%%%%%%%%%%%%%%%%%%%%%%%%%%%%%%%%
\end{Proposition}

 {\sl Proof of Proposition~\ref{prop:zrwq}.}
%%%%%%%%%%%%%%%%%%%%%%%%%%%%%%%%%%%%%%%%
Let~$\alpha:=\{2,3,\dots,n \}$ and~$ \beta := \{1,2, \dots n-1\}$.
%%We begin by determining the minimal~$w$ such that~$A^w(\alpha|\beta)$ is positive.  
%%%
Suppose that~$A^w(\alpha|\beta)>0$ for some~$w>0$. 
Since each copy of~$A$ contains~$s-1$ diagonal
 lines corresponding to~$W_2(\ell^2)\cdots W_{s}(\ell^{s})$ 
(see  Figure~\ref{fig:EB_Z1}),
 the~$n-1$ vertex-disjoint paths
corresponding to~$A^{w}(\alpha|\beta)$
 require at least~$\mu_1$ copies of~$A$. We conclude that
\be\label{eq:pollo}
r_{\ell}(A) \geq \mu_1.
\ee

We now  show that~$A^{\mu_1}(\alpha|\beta)>0$.
We introduce more notation. Consider the network corresponding to~$A^{\mu_1}$.
The symbol~$\searrow$
represents a diagonal arc from one of~$W_2(\ell^2),\dots, W_{s}(\ell^{s})$ in some copy of~$A$, and is accompanied by an explanation of which~$ L_j$ is used. We now describe a set of paths in the network corresponding
to~$A^{\mu_1}$. 
 The first path is
\begin{align*}
P_2: \; 2&\to 2 \to \dots \to 2  \searrow 1 \to 1 \to \dots \to 1,
%%%%%%
\end{align*}
where~$\searrow$ is from~$L_2$ in 
 the first copy of~$W_2$. The second path is
\begin{align*}
P_3: \; 3&\to 3 \to \dots \to 3  \searrow 2 \to 2 \to \dots \to 2,
%%%%%%
\end{align*}
where~$\searrow$ is from~$L_3$ in 
 the first copy of~$W_3$ if~$s>2$ or from~$L_3$ in the second copy of~$W_2$ if~$s=2$. 
Note that this implies that~$P_2$ and~$P_3$ are vertex-disjoint. 
The last path is
\begin{align*}
P_{n}: \; n&\to n \to \dots \to n  \searrow n-1 \to n-1 \to \dots \to n-1,
%%%%%%
\end{align*}
where~$\searrow$ is from~$L_n$ in 
 the first possible copy of some~$W_i$ such that~$P_n$ and~$P_{n-1}$
are vertex-disjoint.   Note that all these paths are vertex-disjoint by construction.

The paths~$P_2,\dots,P_{s}$ use the first copy of~$W_2(\ell^2)\cdots W_{s}(\ell^{s})$.
The paths~$P_{s+1},\dots , P_{2s-1}$ use the second copy of~$W_2(\ell^2)\cdots W_{s}(\ell^{s})$, and so on. This implies that~$\mu_1$ copies of~$A$ are indeed enough to realize this family of~$n-1$ vertex-disjoint paths,
so~$A^{\mu_1}(\alpha|\beta) >0$.

Pick~$c\in\{1,\dots,n-2\}$. 
Let~$\gamma:=\{  n-c+1,\dots,n  \}$ and~$\delta:=\{1,\dots,c \}$.
Then~$A(\gamma|\delta)$ is  a
 lower-left corner minor,  $(\gamma|\delta)\not =(\alpha|\beta)$, and~$A(\gamma|\delta)$ is not 
the determinant of~$A$.
 To complete the proof, we need to show  
  that~$A^{\mu_1}(\gamma|\delta)>0$. 
	To do this, we generate a set of  vertex-disjoint paths connecting every~$ x\in\gamma$ to~$(x-n+c )\in \delta$. 
	We do this by modifying the paths~$P_i$ described above. 
	First, delete the paths~$P_n,\dots,P_{c+2}$, as these are not needed now.
	The path~$P_{c+1}$ contains
	an arc connecting vertex~$c+1$ to~$c$. The structure of the diagonal lines in the  
	network (see Figure~\ref{fig:EB_Z1}) implies that we can extend this arc to a path~$\tilde P_{c+1}$
	connecting~$n$ to~$c$, as every diagonal line emanates from~$n$.
	The path~$P_{c }$ contains
	an arc connecting vertex~$c $ to~$c-1$.
	This  path can be extended similarly 
		to a path~$\tilde P_{c}$ connecting~$n-1$ to~$c-1$ 
			that does not intersect nor
	touch~$\tilde P_{c+1}$. Continuing in this 
	fashion shows that~$A^{\mu_1}(\gamma|\delta)>0$.~\hfill{$\square$}

The   second class of parametrized matrices, denoted~$Z_2(s)$,
 requires more notation.
Note that the multiplier of~$L_i$
in~$W_i(\ell^i) $ is~$  \ell^i_{n-i+1} =\ell_{x_i+n-i}$.
For any~$i=2,\dots,n-1$, we 
    define a subset     of~$W_i(\ell^i)$  by
\begin{align}\label{eq:P_i}
%%%%%%%%%%%%%%%%%%%%%%%%%%%%%%%
\Psi_i(\ell^i):=\{W_i(\ell^i) :\ell_{x_i+n-i}>0, \; \prod_{j=x_i}^{x_i+n-i-1} \ell_j=0   \}.
%%%%%%%%%%%%%%%%%%%%%%%%%%
\end{align}
This describes all  the cases where the multiplier of~$L_i$ is positive,
and    at least one of the  other multipliers  is zero. 
 For example, for~$n=5$, 
\begin{align*}
\Psi_2(\ell^2)&=\{L_2, L_3L_2,L_4L_2,L_5L_2,L_4L_3L_2,L_5L_4L_2,L_5L_3L_2\}, \\
\Psi_3(\ell^3)&=\{L_3,L_4L_3,L_5L_3\}, \\
\Psi_4(\ell^4)&=\{L_4\},
\end{align*}
where it is understood that all the multipliers are positive. 

\begin{Definition}\label{def:Z2mats}
$Z_2(s)$, $s\in D_2:=\{2,\dots,n\} $, is the class of matrices in the form:
\be\label{eq:fdrgw}
 L_2(\ell_2)\cdots L_{s-2}(\ell_{s-2})P(\ell^{s-1}) W_s(\ell^s)\cdots W_n(\ell^n),
\ee
where~$  P(\ell^{s-1})\in\Psi_{s-1}(\ell^{s-1})$, and
all the multipliers are positive. 
%%%%%%%%%%%%%%
\end{Definition}
%%%%%%%%%%%%%%%%%%%%%%%%%%

Thus,
 $Z_2$ includes  all the matrices in the form
\begin{align*}
&W_2 W_3 \dots W_n, \\
&{P_2} W_3\cdots W_n,\\
& L_2P_3 W_4\cdots W_n,\\
&\vdots \\
&L_2\cdots L_{n-2}P_{n-1}W_n,
\end{align*}
 where~$P_i\in\Psi_i$. 

A special case is when~$P(\ell^{s-1})=L_{s-1}(\ell^{s-1})$ and then~\eqref{eq:fdrgw} becomes
\[
 L_2(\ell_2)\cdots L_{s-2}(\ell_{s-2})L_{s-1}(\ell_{s-1})  W_s(\ell^s)\cdots W_n(\ell^n).
\]

The   planar network associated 
with~$Z_2(s)$ represents the product  
%\[
$
L_2\cdots L_{s-2}PW_s\cdots W_n,
%\]
$
with~$P\in\Psi_{s-1}$ (we omit the parameters for the sake of clarity).
 Note that~$P$ must include the matrix~$L_{s-1}$, 
and must not include at least one of the matrices~$L_n,\dots,L_{s}$.
%%%%%%%%%%%%%%%%
This network includes horizontal lines,
 and at least~$n-1$ downward  pointing  diagonal lines.
There are~$n-s+1$ 
   diagonal lines corresponding to~$W_x$, $x=s,\dots,n$.
	Each such line 
connects~$n$ to~$x-1$. There are another~$s-3 $ 
  diagonal lines corresponding to~$L_x$,  $x=2,\dots,s-2$.  
	Each such line connects~$x$ to~$x-1$. 
The matrix~$P$ corresponds to a diagonal arc that 
 connects~$s-1$ to~$s-2$, and may also include other diagonal arcs. 
Figure~\ref{fig:EB_Z2} depicts this network for~$n=6$, the particular case~$P=L_{s-1}$, and
 all~$s\in\{2,\dots,6\}$.

 \begin{figure*}
\begin{center}
\includegraphics[width= 13cm,height=7cm]{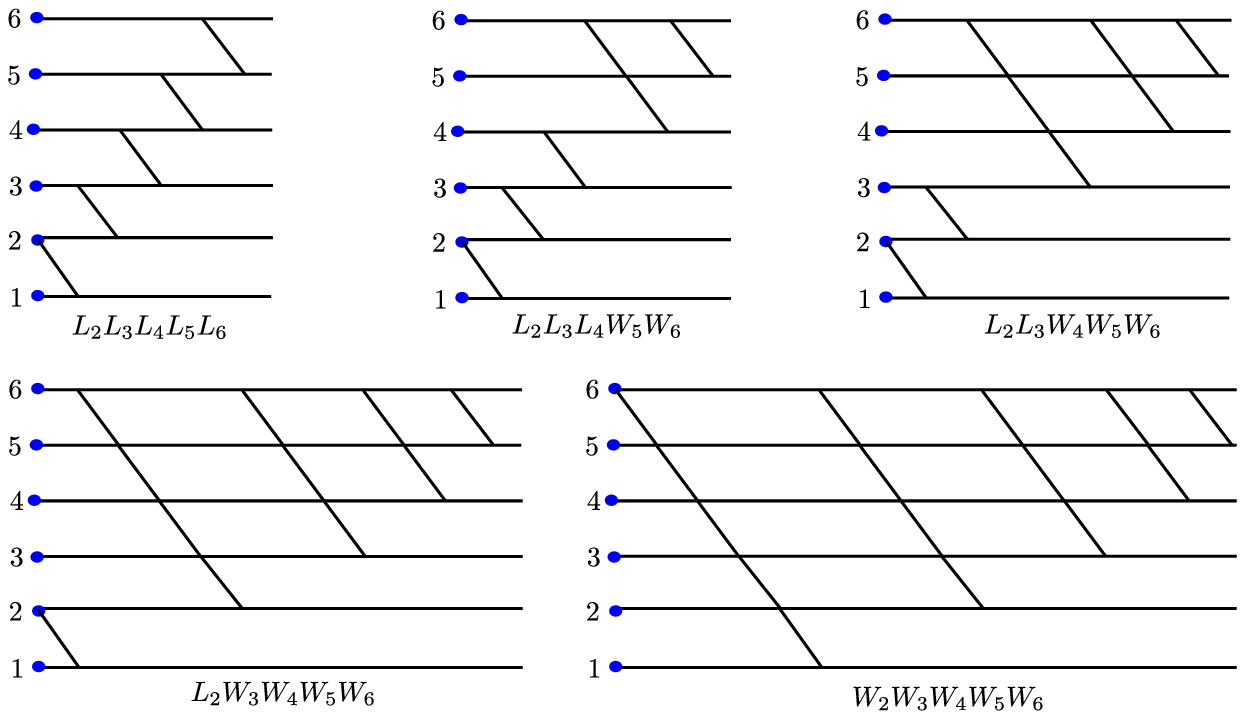}
\caption{\updt{Planar network  corresponding to~$Z_2(s)$  with~$n=6$, the particular case~$P=L_{s-1}$, and~$s\in\{2,\dots,6\}$. Upper-left figure corresponds to~$s=6$, upper-middle to~$s=5$, upper-right to~$s=4$, lower-left to~$s=3$, and lower-right to~$s=2$.}}\label{fig:EB_Z2}
\end{center}
\end{figure*}

\begin{Proposition}\label{prop:zrwq2}
%%%%%%%%%%%%%%%%%%%%%
Pick~$s\in D_2$ and a  matrix~$A\in Z_2(s)$.  Let~$\mu_2:=    s-1   $.
 Then~$r_\ell(A)=\mu_2 $. 
%%%%%%%%%%%%%%%%%%%%%%%%%%%%%%%%%%%%%
\end{Proposition}

{\sl Proof of Proposition~\ref{prop:zrwq2}.}
%%%%%%%%%%%%%%%%%%%%%%%%%%%%%%%%%%%%%%%%%%%%%%%%%%%%%%%%%%%%%%%%%%%%%
We first assume that~$P=L_{n-1}$. 
Let~$\alpha:=\{n\}$ and~$\beta:=\{1\}$.
We claim that  
 \be\label{eq:hfpoi}
\min\{k: A^k(\alpha|\beta)>0 \} = \mu _2 . 
\ee
Indeed,
in the network of (every copy) of~$A$ the longest downward pointing diagonal 
connects~$n$   to~$s-1$. For any vertex in~$\{ s-2,s-3,\dots,1 \}$ (i.e. all the 
vertices ``below''~$s-1$)
the longest  
 downward pointing diagonal is a single arc. 
%%%%%%%%
This implies that at 
least~$ n-(n-s+1) = \mu_2$ copies of~$A$ are needed to realize a path from the source vertex~$n$
to the sink vertex~$1$. Given this number of copies, the path is:
 %%%
\begin{align*}
  n \searrow s-1 \searrow s-2 \searrow s-3 \searrow \cdots \searrow 1. 
\end{align*}
 This proves~\eqref{eq:hfpoi}. 

Pick~$c\in\{2,\dots,n-1\}$. 
Let~$\gamma:=\{  n-c+1,\dots,n \}$ and~$\delta:=\{1,\dots,c \}$.
Then~$(\gamma|\delta)$ is  a
 lower-left corner minor that is not~$(\alpha|\beta)$ nor
the determinant of~$A$. To complete the proof, we need to show  
  that~$A^{\mu_2}(\gamma|\delta)>0$. 
	To do this, we generate a set of  vertex-disjoint paths connecting
	every~$ x\in\gamma$ to~$(x-n+c )\in \delta$ in the network of~$A^{\mu_2}$:
in the first copy, we use
\begin{align*}
										 & n-c+1 \searrow s-1  ,\\
										& n-c +2 \searrow s  ,\\
										&\;\vdots \\
										& n\searrow    c+s-2.
\end{align*}
Note that this is vertex-disjoint set of paths. In the second copy, we use 
\begin{align*}
										 & s-1 \searrow s-2  ,\\
										& s\searrow s-1  ,\\
										&\;\vdots \\
										& c+s-2\searrow    c+s-3.
\end{align*}
and similar ``one-step'' paths in all the next copies. 
Then~$\mu_2=s-1$ copies are always enough for a
 set of  vertex-disjoint paths connecting
	every~$ x\in\gamma$ to~$(x-n+c )\in \delta$, and this completes the proof 
	when~$P=L_{n-1}$.

Now assume that~$P\ne L_{s-1}$, i.e.~$P=QL_{s-1}$, where~$Q$ contains at least one matrix from~$L_n,\dots,L_{s}$, but not all of these matrices. 
In the planar network of~$Q$ the longest possible downward pointing
diagonal   that starts from the vertex~$n$ contains~$n-s$ consecutive edges, 
corresponding to~$Q=L_n\cdots L_{s+1}$. 
Since the downward pointing diagonal line  
 corresponding to~$W_s$ starts from~$n$ and contains~$n-s+1$ consecutive edges, using diagonal lines corresponding to~$Q$ cannot reduce the number of copies of~$A$ required to obtain a path from
 the source vertex~$n$ to the sink vertex~$1$. 
%%% 
This completes the proof.~\hfill{$\square$}.}

The next result follows from Theorem~\ref{thm:deg} and 
Propositions~\ref{prop:zrwq} and~\ref{prop:zrwq2}.

%%%%%%%%%%%%%%%%%%%%%%%%%%%%%%%%%%%%%%%%%%
\begin{Corollary}\label{coro:main_explicit}
% ===================================
Let~$A\in\R^{n\times n}$  be an oscillatory matrix.
If~$A=LDU$ with~$L \in Z_i(s_i)$ for some~$i\in\{1,2\}$,
$D$ a positive diagonal matrix, and~$U^T \in Z_j(s_j)$ for some~$j\in\{1,2\}$,
 then
\be\label{eq:r_explicit}
r(A)=\max\{\mu_i(s_i), \mu_j(s_j)\}.
\ee
\end{Corollary}

Note that~$\mu_1(s)$ takes values in $\{n-1,\left\lceil  \frac{n-1}{2} \right\rceil, \left\lceil  \frac{n-1}{3} \right\rceil,\dots,1 \} $,
and~$\mu_2(s)$    takes values in~$\{1,2 ,\dots, n-1\}$, so
Corollary~\ref{coro:main_explicit}
 covers oscillatory matrices with all possible exponent values.

The following examples demonstrate our theoretical results.
\begin{Example}\label{exp:z1}
%%%%%%%%%%%%%%%%%%%%%%%%%%%%%%%%%%%%%%%%%%%%%%%%%
Suppose that~$L \in Z_1(n)$ and~$U^T\in  Z_1(n)$.
 In this case, all the multipliers in the SEB factorization of~$A=LDU$ are positive, so
  Proposition~\ref{prop:seb_TP} implies that~$A $ is~TP. 
On the other-hand,~\eqref{eq:r_explicit} gives
\[
r(A)=\max\{\mu_1(n), \mu_1(n)\}=\lceil \frac{n-1}{n-1}\rceil=1.
\]

As another example, suppose that~$L\in Z_1(2)$ and~$U^T \in Z_1(2)$.
 Then~$A$ is as in~\eqref{eq:bosc1}, and~\eqref{eq:r_explicit} gives
\[
r(A)=\max\{\mu_1(2), \mu_1(2)\}=\lceil \frac{n-1}{1} \rceil=n-1.
\]
 This recovers one of the results in~\cite{fallat2007class}.
\end{Example}

For a  vector~$x=\begin{bmatrix} x_1&x_2&\dots&x_n \end{bmatrix}^T\in\R^n$,
let~$\underline{x}:=\begin{bmatrix} x_n&x_{n-1}&\dots&x_1\end{bmatrix}^T$.
%%%
\begin{Example}
%%%%%%%%%%%%%%%%%%%%%%%%%%%%%%%%%%%
As a more concrete example, consider
\[
A = 
\begin{bmatrix}
1 & 1 & 2 & 2 \\
2 & 3 & 7 & 9 \\
6 & 9 & 22 & 30 \\
6 & 9 & 22 & 31
\end{bmatrix},
\]
which is I-TN, but not TP (for example,~$A(\{2,3\}|\{1,2\})=0$). The SEB factorization of~$A$ is
\begin{align*}
A&=L_4(1)L_3(3)L_2(2)IU_3(1)U_4(2)U_2(1)U_3(2)U_4(1) \\
&=W_2(\ell^2)Q_3(u^3)Q_2(u^2),
\end{align*}
with~$\ell^2=\begin{bmatrix} 1 & 3 & 2  \end{bmatrix}^T$, $u^3=\begin{bmatrix} 1 & 2 \end{bmatrix}^T$, and~$u^2=\begin{bmatrix} 1 & 2 & 1 \end{bmatrix}^T$, i.e.~$A$ is oscillatory. 
Since~$L=W_2(\ell^2)\in Z_1(2)$ and~$U^T=W_2(\underline{u^2})W_3(\underline{u^3})\in Z_1(3)$,   
\[
r(A)=\max\{\mu_1(2), \mu_1(3)\}=\max\{3,2\}=3.
\]
Indeed, 
\[
A^2=
\begin{bmatrix}
27 & 40 & 97 & 133 \\
104 & 155 & 377 & 520 \\
336 & 501 & 1219 & 1683 \\
342 & 510 & 1241 & 1714
\end{bmatrix},
\] 
and so~$(A^2)(\{2,3,4\}|\{1,2,3\})=0$, thus~$A^2$ is not TP, implying that indeed~$r(A)=3$.
\end{Example}

\begin{Example}
%%%%%%%%%%%%%%%%%%%%%%%%
Let~$L\in Z_2(2)$ and~$U^T\in Z_2(2)$. This also means that~$L\in Z_1(n)$ and~$U^T\in Z_1(n)$, and by Example~\ref{exp:z1}, $r(A)=1$. Indeed Eq.~\eqref{eq:r_explicit} gives 
\[
r(A)=\max\{\mu_2(2), \mu_2(2)\}= 1.
\]
Now assume that~$L\in Z_2(n)$ and~$U^T \in Z_2(n)$. Thus,~$A$ is as in~\eqref{eq:jacobi} (note that~$W_n(\ell^n)=L_n(\ell_k)$ and~$Q_n(u^n)=U_n(u_k)$), and in this case
\[
r(A)=\max\{\mu_2(n), \mu_2(n)\}=n-1\]
 (see~\cite{fallat2007class}).

As a concrete example, consider
\[
A = 
\begin{bmatrix}
1 & 2 & 0 & 0 & 0 \\
2 & 5 & 3 & 0 & 0 \\
0 & 2 & 7 & 2 & 6 \\
0 & 8 & 29 & 11 & 34 \\
0 & 24 & 89 & 41 & 131
\end{bmatrix},
\]
which is I-TN, but clearly not TP. The SEB factorization of~$A$ is 
\begin{align*}
A&=L_2(2)L_5(3)L_4(4)L_3(2)L_5(2)L_4(1)L_5(2)IU_5(1)U_4(2)U_5(3)U_3(3)U_2(2) \\
&=L_2(2)W_3(\ell^3)W_4(\ell^4)W_5(\ell^5)Q_5(u^5)Q_4(u^4)U_3(3)U_2(2),
\end{align*}
with~$\ell^3=\begin{bmatrix} 3 & 4 & 2  \end{bmatrix}^T$, $\ell^4=\begin{bmatrix} 2 & 1 \end{bmatrix}^T$,  $\ell^5=2$, $u^5=1$, and~$u^4=\begin{bmatrix} 2 & 3 \end{bmatrix}^T$, i.e.~$A$ is oscillatory.
Since~$L=L_2(2)W_3(\ell^3)W_4(\ell^4)W_5(\ell^5)\in Z_2(3)$ and~$U^T=L_2(2)L_3(3)W_4(\underline{u^4})W_5(\underline{u^5})\in Z_2(4)$, 
\[
r(A)=\max\{\mu_2(3), \mu_2(4)\}=\max\{2,3\}=3.
\]
Since~$A^2(1,4)=0$, $A^2$ is not~TP. It can be verified that~$A^3$ is TP, implying that indeed~$r(A)=3$.
\end{Example}

We now describe several  generalizations of our results. 
%%%%%%%%%%%%%%%%%%%%%%%%%%%%%%%%%%%%%%%%
\subsection{Generalizations}
%%%%%%%%%%%%%%%%%%%%%%%%%%%%%%%%
\updt{
In many SEB factorizations of oscillatory matrices,~$L$ or~$U^T$ can be \emph{rewritten}
 as a class in~$Z_1$ or~$Z_2$.
We demonstrate this using  an example.  Consider the case~$n=5$ and the product
\[
U:=[U_5(u_1)][U_4(u_2)U_5(u_3)][U_3(u_4)U_4(u_5)U_5(u_6)][U_2(u_7)U_3(u_8)U_4(u_9)U_5(u_{10})],
\]
with~$u_3=u_5=u_8=u_{10}=0$, and all the other multipliers positive, that is,
\[
U=U_5(u_1)U_4(u_2)U_3(u_4)U_5(u_6)U_2(u_7)U_4(u_9).
\]
This seems to suggest that~$U^T \notin (Z_1 \cup Z_2)$. 
However, using~\eqref{eq:Lij_flip} yields
%%%%%%%%%%%%%%%%%%%%%%%%%%%%%%%%%% 
\begin{align*}
U&=U_5(u_1)U_4(u_2)U_5(u_6)U_3(u_4)U_4(u_9)U_2(u_7)\\
&=Q_5(u^5)Q_4(u^4)U_3(u_4)U_4(u_9)U_2(u_7),
\end{align*}
with
\[
u^5:=u_1, \; u^4:=\begin{bmatrix} u_2  &u_6 \end{bmatrix}^T.
\]
Thus,~$U^T\in Z_2(4)$ (with~$P=L_4(u_9)L_3(u_4)$).
}
%%%%%%%%%%%%%%%%%%%%%%%%%%%
 
%%%%%%%%%%%%%%%%%%%%%%%%%%%
\updt{
This suggests  that in many factorizations, $L$ and~$U^T$ can be rewritten as classes in~$Z_1 \cup Z_2 $,
 and thus the corresponding exponent can be determined explicitly.
}

Our analysis uses  the fact that certain edges in the associated planar networks have positive weights, while ignoring their actual value. Thus, another  generalization is to
 arbitrary products of oscillatory matrices.
 
 \updt{
\begin{Corollary}\label{coro:main_prod}
%%%%%%%%%%%%%%%%%%%%%%%%%%%%%%%%%%%%%%%%%%%%%%%%%%%%%
 Consider a  set of matrices~$\{A_x\}_{x=1}^m$, where every~$A_x\in\R^{n \times n}$   is  oscillatory.
 Let~$L(x)$ [$U(x)$] denotes the product of all~$L_i$ [$U_i$] matrices in the SEB factorization of~$A_x$. Suppose that there exist~$i,j,s_i,s_j$ such that~$L(x) \in Z_i(s_i)$ and~$(U(x))^T\in Z_j(s_j)$ for all~$x$, but possibly every~$L(x)$ and every~$U(x)$ has  a different set of multiplier values.
 Then a  product of any~$k$ matrices from the set is~TP
iff
\[
k\ge \max\{\mu_i(s_i),\mu_j(s_j)\}.
\]
 \end{Corollary}
 } % of \updt
 
The next example demonstrates this.  
\begin{Example}
Let
\[
A_1=
\begin{bmatrix}
2 & 8 & 16 & 48 \\
8 & 33 & 67 & 203 \\
8 & 39.5 & 89.5 & 289.5 \\
20 & 134.5 & 350.5 & 1219.5
\end{bmatrix}, \quad 
A_2 = 
\begin{bmatrix}
1 & 2 & 8 & 24 \\
2 & 6 & 28 & 88 \\
6 & 23 & 117 & 376 \\
30 & 145 & 789 & 2580
\end{bmatrix}.
\]
The SEB factorizations   are:
\begin{align*}
A_1&=L_4(2.5)L_3(1)L_2(4)L_4(5.5)L_3(6.5)L_4(1)D_1U_3(1)U_4(2)U_2(4)U_3(2)U_4(3), \\
A_2&= L_4(5)L_3(3)L_2(2)L_4(6)L_3(2.5)L_4(2)D_2U_3(2)U_4(1)U_2(2)U_3(4)U_4(3), 
\end{align*}
with
$D_1:=\diag(2,1,3,5)$ and~$D_2:=\diag(1,2,3,4)$.
 This implies that~$A_1$ and~$A_2$ are oscillatory matrices,
 but not~TP (see Propositions~\ref{prop:seb_TP} and~\ref{prop:seb_OSC}). 
%%%
Note that,
\begin{align*}
L(1)&=L_4(2.5)L_3(1)L_2(4)L_4(5.5)L_3(6.5)L_4(1), \\
L(2)&=L_4(5)L_3(3)L_2(2)L_4(6)L_3(2.5)L_4(2), \\
U(1)&=U_3(1)U_4(2)U_2(4)U_3(2)U_4(3), \\
U(2)&=U_3(2)U_4(1)U_2(2)U_3(4)U_4(3).
\end{align*}
This implies  that both~$L(1)$ and~$L(2)$ [$(U(1))^T$ and~$(U(2))^T$] belong to~$Z_1(4)$ [$Z_1(3)$], 
where each has a different set of parameter values.
 Thus,
\[
r(A_1)=r(A_2)=\max\{\mu_1(4),\mu_1(3)\}=\max\{1,2\}=2.
\]
% of \updt
It is straightforward to verify that indeed~$A_1^2$, $A_2^2$, $A_1A_2$, and~$A_2A_1$ are all TP matrices. \updt{Note that since~$\mu_1(4)=1$, our theoretical results show that
all the lower-left corner minors
 (except perhaps for the determinant) of both~$A_1$ and~$A_2$ are positive.}
\end{Example}

\updt{
Note that in  general given a set~$\{A_x\}_{x=1}^m$ of oscillatory matrices,
with~$r(A_x)=w$ for all~$x=1,\dots,m$, it is possible that a product   
  of less than~$w$ 
matrices from the set is~TP. This is demonstrated in the following example.
} % of \updt

%%%%%%%%%%%%%%%%
\begin{Example}\label{exp:diffslsu}
Consider the matrices
\[
A_1 =
\begin{bmatrix}
1 & 1 & 2 & 6 \\
3 & 4 & 9 & 28 \\
3 & 4 & 10 & 32 \\
6 & 8 & 20 & 65
\end{bmatrix}, \quad
A_2 = 
\begin{bmatrix}
1 & 3 & 3 & 6 \\
2 & 7 & 7 & 14 \\
6 & 23 & 24 & 48 \\
6 & 25 & 27 & 55
\end{bmatrix}.
\] 
Their factorizations are
\begin{align*}
A_1&= L_4(2)L_3(1)L_2(3)U_3(1)U_4(1)U_2(1)U_3(2)U_4(3), \\
A_2 &= L_4(1)L_3(3)L_2(2)L_4(1)L_3(2)U_2(3)U_3(1)U_4(2).
\end{align*}
\updt{
Thus,~$L(1)\in Z_1(2)$, $(U(1))^T\in Z_1(3)$, $L(2)\in Z_1(3)$, and~$(U(2))^T\in Z_1(2)$.Thus,
\begin{align*}
r(A_1)=\max\{\left\lceil  \frac{4-1}{2-1} \right\rceil,\left\lceil  \frac{4-1}{3-1} \right\rceil\}=3,
\end{align*}
 and
\[
r(A_2)=\max\{\left\lceil  \frac{4-1}{3-1} \right\rceil,\left\lceil  \frac{4-1}{2-1} \right\rceil\}=3.
\]
However, it can be verified that both~$A_1A_2$ and~$A_2A_1$ are~TP, i.e.  we do not necessarily 
need a product of three matrices  
to obtain a~TP matrix.
} % of \updt
%%%%
\end{Example}

\subsection{Upper-Bounds on~$r(A)$}
%%%%%%%%%%%%%%%%%%%%%%%%%%
\updt{
It is clear that adding new edges to a planar network associated with an oscillatory matrix can never increase its exponent value. Thus, if a given class of oscillatory matrices can be factored 
as matrices in~$Z_i$ multiplied by additional terms, 
then our results can be used to derive an upper bound on
their exponent. This is stated formally in the next result. 
   
%%%%%%%%%%%%%%%%%%%%%%%%%%%%%%%%%%%%
\begin{Corollary}\label{corr:bound}
%%%%%%%%%%%%%%%%%%%%%%%%%%%%%%%%%%%%
Let~$A=LDU\in\R^{n\times n}$ be an oscillatory matrix. 
If there exist~$i,j$ such that either~$L \in Z_i(s_i)$ and~$U^T$ 
belongs  to~$Z_j(s_j)$ multiplied by  some additional terms, 
or~$U^T \in Z_j(s_j)$ and~$L$
belongs  to~$Z_i(s_i)$ multiplied by  some additional terms,
 then
\[
r(A)\le \max\{\mu_i(s_i),\mu_j(s_j)\}.
\]
\end{Corollary}
%%%
 For example, suppose that~$A\in\R^{5\times 5} $ can be  factored as 
\[
A=L_3(\ell_3)L_2(\ell_4)L_3(\ell_7)W_4(\ell^4)W_5(\ell^5)DQ_5(u^5)Q_4(u^4)U_3(u_7)U_2(u_4).
\]
Then~$L=L_3(\ell_3)W$, where~$W:=L_2(\ell_4)L_3(\ell_7)W_4(\ell^4)W_5(\ell^5)$, and $U^T \in L_2(u_4)L_3(u_7)W_4(u^4)W_5(u^5)$. Here~$U^T \in Z_2(4)$,
 and since~$W \in Z_2(4)$, we see that~$L$ belongs  to~$Z_2(4)$ multiplied by 
the additional term~$L_3(\ell_3)$. Thus,~$r(A) \le \mu_2(4)=3$.
}

The mapping from the SEB factorization of~$A$ to~$r(A)$ is nontrivial. 
There are known cases  where adding specific
additional terms necessarily reduces the  exponent~\cite{fallat2007class}.
\updt{
%%%%
But in general  adding terms with positive multipliers
to~$L\in Z$ or to~$U^T\in Z$ does not necessarily decrease~$r(A)$
 and thus in these cases the bound in Corollary~\ref{corr:bound}
may be tight. For example, if~$r(A)=2$
then
adding terms to the SEB factorization of~$A$ 
 cannot decrease~$r(A)$ unless the modified factorization includes all the multipliers in~\eqref{eq:EB}, as the SEB factorization of a TP matrix is unique. The following example demonstrates this. 
%The following example demonstrates that in general adding terms with positive multipliers to~$L\in Z$ or to~$U^T\in Z$ does not necessarily reduce~$r(A)$ and thus in these cases the bound in Corollary~\ref{corr:bound} may be tight (see also Remark~\ref{rmk:mu4_s_even}).  
%%%%%%%%%%%%%%%%%%%%%%%%%%%%%%%%%%
\begin{Example}
Let
\[
A = 
\begin{bmatrix}
1 & 1 & 4 \\
1 & 2 & 10 \\
0 & 2 & 13
\end{bmatrix}.
\]
Its SEB factorization is~$L_2(1)L_3(2)IU_3(2)U_2(1)U_3(4)$, i.e.~$L \in Z_2(3)$ and~$U^T$ belongs 
to~$Z_2(3)$ multiplied by  the additional term~$L_3(4)$. 
By Corollary~\ref{corr:bound}, 
\[r(A)\leq \mu_2(3)=2 ,
\]
 and indeed in this case this bound is tight.
\end{Example}

} % \of \updt
\section{Conclusion}\label{sec:conc}
 %%%%%%%%%%%%%%%%%%%%%%%%
%Oscillatory matrices are a class of   matrices that are intermediary between TN and TP matrices. 
The exponent~$r(A) $ of an 
oscillatory matrix~$A\in\R^{n\times n}$ is the least positive integer such that~$A^r$ is~TP.  
It is well-known that~$r(A)\in\{1,\dots,n-1\}$.
\updt{Fallat and Liu}~\cite{fallat2007class} used  the SEB factorization,  and its associated planar network,
to describe classes  of oscillatory matrices satisfying~$r(A)=n-1$.

\updt{
Here, we use the planar network
to derive
  explicit expressions for~$r(A)$ for several classes of oscillatory matrices \updt{with exponents between~$1$ and~$n-1$}. In addition, we provide non-trivial upper-bounds on~$r(A)$ for several other classes of oscillatory matrices.  
	
	Our analysis is based on writing~$r(A)$ as the maximum of~$r_\ell(A)$ and~$r_u(A)$
	and finding matrices for which these two terms can be determined explicitly using the SEB factorization.  
	The latter is based on a ``worst-case'' approach, that is, 
	finding a specific corner minor that requires the maximal number of copies of~$A$ to make it positive. 
	
		  We believe that the SEB factorization and the associated planar network are powerful tools for analyzing I-TN matrices in general, and will find more applications in the 
analysis of  oscillatory matrices.
%%%%
Possible topics for further research include the following. First, a natural extension to this work is identifying additional  classes of oscillatory matrices whose exponent can be determined explicitly.  Specifically, given the SEB factorization of any oscillatory matrix, can its exponent be determined explicitly? 

Also, Corollary~\ref{coro:main_prod} considers   arbitrary products of oscillatory matrices that share the same SEB factorization, possibly with different multiplier values and different diagonal matrices.  Can this result  be extended to product of matrices that share only their exponent value? As suggested by Example~\ref{exp:diffslsu}, these cases seem to be different than the result in Corollary~\ref{coro:main_prod}.
} % of \updt
Finally, let~$\lambda_1>\lambda_2>\cdots>\lambda_n$ denote the eigenvalues of~$A\in\R^{n\times n}$. 
Then for large~$m$ we have~$  \Trace(A^m) =   \lambda_1^m+\cdots +\lambda_n^m  \approx
\lambda_1^m$, where~$\Trace(B)$ denotes the sum of the diagonal entries of~$B$. Assume that~$A$ is I-TN, and thus it admits an SEB factorization. Note that~$A^m$ also admits an SEB factorization (and its associated planar network in obtained by concatenating~$m$ times the planar network associated with~$A$). Using the SEB factorization of~$A^m$ (and the associated planar network), what can be said 
 about~$\Trace(A^m)$?

 %%%%%%%%%%%%%%%%%%%%%%%%%%%%%%%%%%%
 \section*{Declaration of competing interest}
 %%%%%%%%%%%%%%%%%%%%%%%%
 There are no competing interests.

  %%%%%%%%%%%%%%%%%%%%%%%%%%%%%%%%%%%
 \section*{Acknowledgments}
 %%%%%%%%%%%%%%%%%%%%%%%%
 We are very grateful to one of the  anonymous referees for many
    comments that greatly helped us in improving this paper.

%%%%%%%%%%%%%%%%%%%%%%%%%%%%%%%%%%%
%\bibliography{Osci_bib,lolo_yz}
%%%%%%%%%%%%%%%%%%%%%%%%%%%%%%%%%%%

 \end{document}